\documentclass[12pt,a4paper,twoside]{article}

\newcommand{\pageformat}[6]{\setlength{\hoffset}{-1in}
                  \setlength{\voffset}{-1in}
                  \addtolength{\hoffset}{#5}
                            \addtolength{\voffset}{#6}
                            \setlength{\oddsidemargin}{#1}
                            \setlength{\evensidemargin}{#2}
                            \setlength{\textwidth}{\paperwidth}
                  \addtolength{\textwidth}{-\oddsidemargin}
                  \addtolength{\textwidth}{-\evensidemargin}
                  \addtolength{\textwidth}{-\marginparsep}
                  \addtolength{\textwidth}{-\marginparwidth}
                            \setlength{\topmargin}{#3}
                            \setlength{\textheight}{\paperheight}
                  \addtolength{\textheight}{-\topmargin}
                  \addtolength{\textheight}{-\headheight}
                  \addtolength{\textheight}{-\headsep}
                  \addtolength{\textheight}{-\footskip}
                  \addtolength{\textheight}{-#4}}
\pageformat{2cm}{3cm}{25mm}{25mm}{1pt}{0pt}

\usepackage{ifthen}
\newboolean{article}
    \setboolean{article}{true}
\newboolean{report}
\newboolean{book}
\newboolean{letter}
\newboolean{italian}
\newboolean{german}
\newboolean{nobaselinestretch}
\newboolean{nosectionappendix}
\newboolean{oldtoc}
\newboolean{nosectionequn}
    \setboolean{nosectionequn}{true}
\newboolean{notheorem}

\ifthenelse{\boolean{german}}{
    \usepackage{german}}{}

\usepackage[latin1]{inputenc}
\usepackage{hyperref}

\usepackage{amsmath}
\usepackage{amssymb}
\usepackage[mathscr]{eucal}

\usepackage{txfonts}

\ifthenelse{\boolean{notheorem}}{}{
    \usepackage{theorem}}

\ifthenelse{\boolean{nobaselinestretch}}{}{
    \renewcommand{\baselinestretch}{1.25}}

\newenvironment{env}[2]{\begin{#1}#2\end{#1}}{}
    \newcommand{\beq}[1]{\begin{env}{equation}{#1}}
    \newcommand{\beqn}[1]{\begin{env}{equation*}{#1}}
    \newcommand{\bal}[1]{\begin{env}{align}{#1}}
    \newcommand{\baln}[1]{\begin{env}{align*}{#1}}
    \newcommand{\bga}[1]{\begin{env}{gather}{#1}}
    \newcommand{\bgan}[1]{\begin{env}{gather*}{#1}}
    \newcommand{\bflal}[1]{\begin{env}{flalign}{#1}}
    \newcommand{\bflaln}[1]{\begin{env}{flalign*}{#1}}
    \newcommand{\bmu}[1]{\begin{env}{multline}{#1}}
    \newcommand{\bmun}[1]{\begin{env}{multline*}{#1}}
    \newcommand{\bsp}[1]{\begin{env}{split}{#1}}

    \newcommand{\eeq}{\end{env}}
    \newcommand{\eeqn}{\end{env}}
    \newcommand{\eal}{\end{env}}
    \newcommand{\ealn}{\end{env}}
    \newcommand{\ega}{\end{env}}
    \newcommand{\egan}{\end{env}}
    \newcommand{\eflal}{\end{env}}
    \newcommand{\eflaln}{\end{env}}
    \newcommand{\emu}{\end{env}}
    \newcommand{\emun}{\end{env}}
    \newcommand{\esp}{\end{env}}

\newcommand{\lf}{\vspace{2ex}}

\renewcommand{\bf}[1]{\textbf{#1}}
\renewcommand{\it}[1]{\textit{#1}}

\renewcommand{\sf}[1]{\textsf{#1}}

\renewcommand{\tt}[1]{\texttt{#1}}
\newcommand{\hl}[1]{\bf{\it{#1}}}
\newcommand{\mrm}[1]{\mathrm{#1}}
\newcommand{\mbf}[1]{\mathbf{#1}}

\newcommand{\cmc}[1]{\mathcal{#1}}
\newcommand{\eus}[1]{\mathscr{#1}}
\newcommand{\euf}[1]{\mathfrak{#1}}
\newcommand{\bb}[1]{\mathbb{#1}}
\newcommand{\nbd}[1]{$#1$\nobreakdash--}
\newcommand{\ol}[1]{\overline{#1}}

\newcommand{\vt}{\vartheta}

\newcommand{\om}{\omega}

\newcommand{\bfam}[1]{\bigl(#1\bigr)}
\newcommand{\Bfam}[1]{\Bigl(#1\Bigr)}
\newcommand{\AB}[1]{\langle#1\rangle}

\newcommand{\CB}[1]{\{#1\}}
\newcommand{\bCB}[1]{\bigl\{#1\bigr\}}
\newcommand{\BCB}[1]{\Bigl\{#1\Bigr\}}

\newcommand{\RO}[1]{[#1)}

\newcommand{\sbars}[1]{\:\bar{#1}^s\:}
\newcommand{\sodots}{\sbars{\odot}}
\newcommand{\set}[2][]{
    \ifthenelse{\equal{#1}{}}{
        \CB{#2}}{
        \CB{#1~|~#2}}}
\newcommand{\bset}[2][]{
    \ifthenelse{\equal{#1}{}}{
        \bCB{#2}}{
        \bCB{#1~|~#2}}}
\newcommand{\Bset}[2][]{
    \ifthenelse{\equal{#1}{}}{
        \BCB{#2}}{
        \BCB{#1~\big|~#2}}}

\DeclareMathOperator{\ls}{\normalfont\sf{span}}
\DeclareMathOperator{\cls}{\ol{\ls}}

\DeclareMathOperator{\id}{\normalfont\sf{id}}

\newcommand{\C}{\bb{C}}

\newcommand{\R}{\bb{R}}

\newcommand{\cA}{\cmc{A}}
\newcommand{\cB}{\cmc{B}}
\newcommand{\cC}{\cmc{C}}

\newcommand{\sB}{\eus{B}}

\newcommand{\sF}{\eus{F}}

\newcommand{\sK}{\eus{K}}

\newcommand{\eH}{\euf{H}}

\newcommand{\U}{\mbf{1}}
\newcommand{\F}{{\mrm{F}}}

\ifthenelse{\boolean{nosectionequn}}{}{
    \numberwithin{equation}{section}
    }

\ifthenelse{\boolean{article}\or\boolean{letter}\or\boolean{nosectionequn}}{
    \setboolean{nosectionappendix}{true}}{}
\ifthenelse{\boolean{nosectionappendix}}{}{
    \renewcommand{\appendix}{
        \chapter*{\appendixname}
        \addcontentsline{toc}{chapter}{\appendixname}
        \renewcommand{\thesection}{\Alph{section}}
        \setcounter{section}{0}}}
   
\ifthenelse{\boolean{report}\or\boolean{book}}{
    }{}

\ifthenelse{\boolean{notheorem}}{}{
        \newcommand{\definame}{Definition.}
        \newcommand{\propname}{Proposition.}
        \newcommand{\lemname}{Lemma.}
        \newcommand{\exname}{Example.}
        \newcommand{\exername}{Exercise.}
        \newcommand{\remname}{Remark.}
        \newcommand{\obname}{Observation.}
        \newcommand{\thmname}{Theorem.}
        \newcommand{\corname}{Corollary.}
        \newcommand{\proofname}{Proof.}
    \ifthenelse{\boolean{german}}{
        \renewcommand{\exname}{Beispiel.}
        \renewcommand{\exername}{Übung.}
        \renewcommand{\remname}{Bemerkung.}
        \renewcommand{\obname}{Beobachtung.}
        \renewcommand{\thmname}{Satz.}
        \renewcommand{\corname}{Korollar.}
        \renewcommand{\proofname}{Beweis.}}{}
    \ifthenelse{\boolean{italian}}{
        \renewcommand{\definame}{Definizione.}
        \renewcommand{\propname}{Proposizione.}
        \renewcommand{\exname}{Esempio.}
        \renewcommand{\exername}{Esercizio.}
        \renewcommand{\remname}{Nota.}
        \renewcommand{\obname}{Osservazione.}
        \renewcommand{\thmname}{Teorema.}
        \renewcommand{\corname}{Corollario.}
        \renewcommand{\proofname}{Dimostrazione.}

       \renewcommand{\appendixname}{Appendice}

       }{}
    \theoremheaderfont{\normalfont\bfseries}
    \theoremstyle{change}
        \theorembodyfont{\rmfamily}
            \newtheorem{emp}{}[section]
                \newcommand{\bemp}[1][]{
                    \begin{emp}\hskip-\labelsep\bf{#1}\hskip\labelsep}
                \newcommand{\eemp}{\end{emp}}
\newtheorem{itemp}[emp]{}
                \newcommand{\bitemp}[1][]{
                    \begin{itemp}\hskip-\labelsep\bf{#1}\hskip\labelsep\normalfont\itshape}
                \newcommand{\eitemp}{\end{itemp}}
            \newtheorem{ex}[emp]{\exname}
                \newcommand{\bex}{\begin{ex}}
                \newcommand{\eex}{\end{ex}}
            \newtheorem{exer}[emp]{\exername}
                \newcommand{\bexer}{\begin{exer}}
                \newcommand{\eexer}{\end{exer}}
            \newtheorem{defi}[emp]{\definame}
                \newcommand{\bdefi}{\begin{defi}}
                \newcommand{\edefi}{\end{defi}}
            \newtheorem{rem}[emp]{\remname}
                \newcommand{\brem}{\begin{rem}}
                \newcommand{\erem}{\end{rem}}
            \newtheorem{ob}[emp]{\obname}
                \newcommand{\bob}{\begin{ob}}
                \newcommand{\eob}{\end{ob}}
        \theorembodyfont{\normalfont\itshape}
            \newtheorem{thm}[emp]{\thmname}
                \newcommand{\bthm}{\begin{thm}}
                \newcommand{\ethm}{\end{thm}}
            \newtheorem{prop}[emp]{\propname}
                \newcommand{\bprop}{\begin{prop}}
                \newcommand{\eprop}{\end{prop}}
            \newtheorem{cor}[emp]{\corname}
                \newcommand{\bcor}{\begin{cor}}
                \newcommand{\ecor}{\end{cor}}
            \newtheorem{lem}[emp]{\lemname}
                \newcommand{\blem}{\begin{lem}}
                \newcommand{\elem}{\end{lem}}
    \newcommand{\qedsymbol}{~\rule[-0.35mm]{2mm}{2mm}}
    \newcounter{proof}[emp]
    \newenvironment{Proof}[1]{
        \vspace{1ex}
        \renewcommand{\item}[1][\stepcounter{proof}(\roman{proof})]%
            {##1\hskip\labelsep}
        \noindent\textsc{#1\hskip\labelsep}}{
        \nolinebreak\qedsymbol}
    \newcommand{\proof}[1][\proofname]{
        \begin{Proof}{#1}\ignorespaces}
    \newcommand{\qed}{\end{Proof}}
    \newcommand{\noqed}{
        \renewcommand{\qedsymbol}{}
        \end{Proof}}}
    \ifthenelse{\boolean{italian}}{
        \renewcommand{\proofname}{Dimostrazione.}}{}

\usepackage[matrix,arrow]{xy}

\begin{document}

\title{Three Ways to Representations of $\sB^a(E)$}
\author{}
\author{
~\\
Michael Skeide\thanks{This work is supported by research fonds of the Department S.E.G.e S. of the University of Molise and by the University of Iowa.}\\\\
{\small\itshape Università degli Studi del Molise}\\
{\small\itshape Dipartimento S.E.G.e S.}\\
{\small\itshape Via de Sanctis}\\
{\small\itshape 86100 Campobasso, Italy}\\
{\small{\itshape E-mail: \tt{skeide@math.tu-cottbus.de}}}\\
{\small{\itshape Homepage: \tt{http://www.math.tu-cottbus.de/INSTITUT/lswas/\_skeide.html}}}\\
}
\date{December 2003}

{
\renewcommand{\baselinestretch}{1}
\maketitle

\begin{abstract}
We describe three methods to determine the structure of (sufficiently continuous) representations of the algebra $\sB^a(E)$ of all adjointable operators on a Hilbert \nbd{\cB}module $E$ by operators on a Hilbert \nbd{\cC}module. While the last and latest proof is simple and direct and new even for normal representations of $\sB(H)$ ($H$ some Hilbert space), the other ones are direct generalizations of the representation theory of $\sB(H)$ (based on Arveson's and on Bhat's approaches to product systems of Hilbert spaces) and depend on technical conditions (for instance, existence of a unit vector or restriction to von Neumann algebras and von Neumann modules). We explain why for certain problems the more specific information available in the older approaches is more useful for the solution of the problem.
\end{abstract}

}



\section{Introduction}\label{intro}

A normal unital representation $\vt\colon\sB(H)\rightarrow\sB(K)$ of the algebra $\sB(H)$ of all adjointable (and, therefore, bounded and linear) operators on a Hilbert space $H$ by operators on a Hilbert space $K$ factors $K$ into the tensor product of $H$ and another Hilbert space $\eH$, such that elements $a$ in $\sB(H)$ act on this tensor product in the natural way by ampli(-fic-)ation, i.e.\
\baln{
K
&
~=~
\eH\otimes H
&&
\text{or}
&
K
&
~=~
H\otimes\eH
\intertext{with}
\vt(a)
&
~=~
\id_\eH\otimes a
&&
\text{or}
&
\vt(a)
&~=~
a\otimes\id_\eH.
}\ealn
It is the goal of these notes to report three different proofs of the following analogue result for Hilbert modules and dicuss their interrelations.

\bitemp[Theorem (Muhly, Skeide and Solel \cite{MSS03p}).]\label{mainthm}
Let $E$ be a Hilbert module over a $C^*$--al\-ge\-bra $\cB$ and let $F$ be a Hilbert module over a \nbd{C^*}algebra $\cC$. If $\vt\colon\sB^a(E)\rightarrow\sB^a(F)$ is a unital strict homomorphism, then there exists a Hilbert \nbd{\cB}\nbd{\cC}module $F_\vt$ and a unitary $u\colon E\odot F_\vt\rightarrow F$ such that $\vt(a)=u(a\odot\id_{\F_\vt})u^*$.

The same result is true, if we replace \nbd{C^*}algebras by \nbd{W^*}algebras, Hilbert modules by \nbd{W^*}modules (and their tensor products) and $\vt$ by a normal homomorphism.
\eitemp

Here, the \hl{strict topology} of $\sB^a(E)$ is the strict topology inherited by considering $\sB^a(E)$ as multplier algebra of the \nbd{C^*}algebra $\sK(E)$ of \hl{compact operators}, which is the norm completion of the \nbd{*}algebra $\sF(E)$ of \hl{finite rank operators} spanned by the \hl{rank-one operators} $xy^*\colon z\rightarrow x\AB{y,z}$. A linear mapping is \hl{strict} (and, therefore, bounded), if it is strictly continuous on bounded subsets of $\sB^a(E)$.

The proof from \cite{MSS03p}, being both the simplest available and the most general, is based on the observation that the tensor product $E\odot E^*$ of the \nbd{\sK(E)}\nbd{\cB}module $E$ and the \hl{dual} \nbd{\cB}\nbd{\sK(E)}module $E^*$ (with inner product $\AB{x^*,y^*}=xy^*\in\sK(E)$ and module operations $bx^*a=(a^*xb^*)^*$) may be identified with $\sK(E)$. (The canonical identification is $x\odot y^*\mapsto xy^*$.) Therefore, since $\vt$ is strict and since $\sK(E)$ has a bounded approximate unit (converging strictly to $\id_E$), we have
\beqn{
F
~=~
\sK(E)\odot F
~=~
(E\odot E^*)\odot F
~=~
E\odot(E^*\odot F)
~=~
E\odot F_\vt
}\eeqn
where we set $F_\vt:=E^*\odot F$. The canonical identification is
\beqn{
E\odot F_\vt
~=~
E\odot(E^*\odot F)
~\ni~
x\odot(y^*\odot z)
~\longmapsto~
\vt(xy^*)z
~\in~
F.
}\eeqn
Clearly, $\vt(a)=a\odot\id_{F_\vt}$.

\brem
A more detailed version can be found in \cite{MSS03p}. The mechanism of the proof can be summarized by observing that, if $E$ is \hl{full} (i.e.\ if the range of the inner product of $E$ generates $\cB$ as a \nbd{C^*}algebra), then $E$ may be viewed as \it{Morita equivalence} from $\sK(E)$ to $\cB$. (If $E$ is not full, then replace $\cB$ by the closed ideal in $\cB_E$ in $\cB$ generated by the inner product.) Then $\sK(E)=E\odot E^*$ and $\cB=E^*\odot E$ serve as identities under tensor product of bimodules. The identifications of the bimodule $F_\vt$ and of $E\odot F_\vt$ with $F$ are highly unique. For instance, we may establish the equality $F=E\odot E^*\odot F$ by showing that $F$ furnished with the embedding $i\colon E\times E^*\times F\rightarrow F$, $i(x,y,z)=\vt(xy^*)z$ has the \it{universal property} of the threefold tensor product $E\odot E^*\odot F$. By these and similar considerations one may see that all identifications are essentially unique by \it{canonical isomorphisms}. We investigate these and other more categorical problems in \cite{MSS03p}. Among the applications of Theorem \ref{mainthm} there is the answer to the question when $\vt$ is a (bistrict) isomorphism, namely, if and only if $F_\vt$ is a Morita equivalence. We will investigate consequences of this insight in Muhly, Skeide and Solel \cite{MSS03p1}.
\erem

Even in the case of normal representations of $\sB(H)$ on another Hilbert space the preceding proof (or, more acurately, its modification to normal mappings) seems to be new. In the remainder, we discuss two known ways of treating the representation theory of $\sB(H)$ (Section \ref{B(G)}). Then we describe modifications to adapt them to Hilbert modules, at least, under certain additional conditions (Sections \ref{Bha} and \ref{Arv}). The two approaches correspond to the two basic constructions of product systems of Hilbert spaces from \nbd{E_0}semigroups on $\sB(H)$, the original one by Arveson \cite{Arv89} based on intertwiner spaces and an alternative one by Bhat \cite{Bha96} based on rank-one operators.

In the generalization to Hilbert modules it turns out that the two product systems constructed by Arveson and by Bhat are well distinguished. The product system constructed by Arveson is, actually, a product system of Hilbert \nbd{\C'}\nbd{\C'}modules where $\C'$ is the commutant of $\C$ when represented in the only possible (non-trivial) way by operators on the Hilbert space $\C$. In terms of the \it{commutant} of Hilbert bimodules (as introduced in Skeide \cite{Ske03c} and also, independently, in Muhly and Solel \cite{MuSo03p}) the Arveson system of an \nbd{E_0}semigroup is the commutant of its Bhat system and the Bhat system is that which corresponds to the representation theory (applied to endomorphisms of $\sB(H)$) in Theorem \ref{mainthm}.

All three proofs of Theorem \ref{mainthm} lead to the construction of product systems of Hilbert (bi-) modules when applied to the endomorphisms of \nbd{E_0}semigroups. We compare the three possibilities. In particular, we emphazise those aspects where the more concrete identifications in Sections \ref{Bha} and \ref{Arv} help solving problems which are more difficult in the above approach. A detailed discussion with complete proofs and specifications about how to distinguish identifications via canonical isomorphism from identifications just via isomorphism can be found in Skeide \cite{Ske03p1}.

\section{Respresentations of $\sB(H)$}\label{B(G)}

In this section we repeat two different ways to look at the representation theory of $\sB(H)$. The goal of this repetition is two-fold. Firstly, it prepairs the terrain for the more subtle arguments in the Hilbert module case. Secondly, we use this opportunity to point at the crucial differences between the two proofs already in the case Hilbert spaces. We hope that the present section will help the reader to understand why these two approaches, whose results may easilly be confused and mixed up in the case of Hilbert spaces, later on, lead to well distinguished directions in the case of von Neumann modules.

Let $H$ denote a Hilbert space and let $\vt$ be a normal unital representation of $\sB(H)$ on another Hilbert space $K$. There are many ways to prove the well-known representation theorem which asserts that there is a Hilbert space $\eH$ such that
\beqn{
K
~\cong~
\eH\otimes H
~\cong~
H\otimes\eH
\text{~~~~~~and~~~~~~}
\vt(a)
~=~
\id_\eH\otimes a
~=~
a\otimes\id_\eH
}\eeqn
in the respective identifications. Which order, $\eH\otimes H$ or $H\otimes\eH$, is the natural one depends heavily on the proof, and the apparent equality of $\eH\otimes H$ or $H\otimes\eH$ is, sometimes, able to cause a certain confusion about the choice.

Following what Arveson \cite{Arv89} did for endomorphisms of $\sB(H)$, we introduce the space of intertwiners
\beqn{
\eH^A=\bCB{x\in\sB(H,K)\colon \vt(a)x=xa~(a\in\sB(H))}.
}\eeqn
One easily checks that $x^*y$ is an element in $\C\U$, the commutant of $\sB(H)$, so that $\AB{x,y}\U=x^*y$ defines an inner product. (Observe that there are well distinguished commutants of $\C$ in $\sB(H)$ and of $\C$ in $\C$.) Of course, being obviously complete, $\eH^A$ is a Hilbert space. A well-known result (for instance, \cite[Lemma 2.10]{MuSo02}) asserts that intertwiner spaces of normal representations act totally, if one of the representations is faithful:
\beq{\label{HAtot}
\cls\eH^AH
~=~
K.
}\eeq
From
\beqn{
\AB{x\otimes h,x'\otimes h'}
~=~
\AB{h,\AB{x,x'}h'}
~=~
\AB{xh,x'h'}
}\eeqn
it follows that $x\otimes h\mapsto xh$ is a unitary $\eH^A\otimes H\rightarrow K$ and that $\vt(a)(xh)=x(ah)$ is the image of $x\otimes ah=(\id_{\eH^A}\otimes a)(x\otimes h)$.

\brem\label{tpiprem}
The reader might find it strange that in the middle term we write $\AB{h,\AB{x,x'}h'}$ instead of $\AB{x,x'}\AB{h,h'}$. However, recalling that $\AB{x,x'}$, acutally is an operator in $\sB(H)'\subset\sB(H)$, the way we wrote it appears, indeed, more natural. Additionally, in the module case only this way of writing remains meaningful and we must dispense with the attitude to put the ``scalars'' outside of the inner product. 
\erem

Although the only von Neumann algebra involved is $\sB(H)$, the preceding proof uses elements from the theory of general von Neumann algebras like the commutant of all operators $a\oplus\theta(a)$ in $\sB(H\oplus K)$ and the fact that bijective algebraic homomorphisms are isomorphisms. Most other proofs make more or less direct use of the fact that a \it{normal} mapping on a von Neumann algebra is known, when it is known on the \it{finite-rank} operators $\sF(H)$ (i.e.\ the subalgebra of $\sB(H)$ spanned by the \it{rank-one} operators $h_1h_2^*\colon h\mapsto h_1\AB{h_2,h}$). One of the most elegant ways to do this we borrow from Bhat \cite{Bha96}. Choosing a reference \it{unit vector} $\om\in H$, we denote by $\eH^B$ the subspace $\vt(\om\om^*)$ of $K$. From
\beq{\label{Hiso}
\AB{\vt(h\om^*)x,\vt(h'\om^*)x'}
~=~
\AB{x,\vt(\om\AB{h,h'}\om^*)x'}
~=~
\AB{h\otimes x,h'\otimes x'}
}\eeq
we see that $h\otimes x\mapsto\vt(h\om^*)x$ defines an isometry $H\otimes\eH^B\rightarrow K$. To see surjectivity we have to make use of an approximate unit for $\sF(H)$ which converges strongly to $\U$ (cf.\ the proof of  Theorem \ref{mainthm} in Section \ref{intro}). Also here we see that $\vt(a)\vt(h\om^*)x=\vt((ah)\om^*)x$ is the image of $ah\otimes x=(a\otimes\id_{\eH^B})(h\otimes x)$.

\brem\label{comprem}
By the uniqueness results mentioned after Theorem \ref{mainthm} the Hilbert space $\eH$ such that $K=H\otimes\eH$ and $\vt(a)=a\otimes\id_\eH$ is unique up to (unique) canonical isomorphism. For instance, the construction of $\eH^B$ depends on the choice of $\om$, but, if $\om'$ is another unit vector, then $\vt(\om'\om^*)\upharpoonright\eH^B$ defines the unique unitary onto the space ${\eH'}^B$ constructed from $\om'$. Moreover, if $\eH=H^*\odot K$ is the Hilbert space according to Theorem \ref{mainthm} (with inner product $\AB{h_1^*\odot k_1,h_2^*\odot k_2}=\AB{k_1,\vt(h_1h_2^*)k_2}$), then $h^*\odot k\mapsto\vt(\om h^*)k$ is the unique unitary $\eH\rightarrow\eH^B$.

We see that $\eH_B$ and $\eH$ are very similar. Indeed, we may say that the construction of $\eH$ is just freeing the construction of $\eH^B$ from the obligation to choose a unit vector.
\erem

\brem
Of course, $\eH^B\cong\eH^A$, but this is an accidental artifact of the fact that $\C'=\C$ and that $\eH^A\otimes H\cong H\otimes\eH^A$, canonically. Considering $\eH^A$ as the \nbd{\C'}\nbd{\C'}module it is, both the expressions $H\otimes\eH$ and $\eH\otimes H$ do not even make sense without additional effort. Indeed, as indicated in Observation \ref{EE'}, to deal with such expressions we have to introduce the tensor product of a Hilbert module over $\cB$ and a Hilbert module over the commutant $\cB'$ of $\cB$.
\erem

\brem
Suppose that $\vt_1$ and $\vt_2$ are unital normal representations of $\sB(H)$ on $K$ and of $\sB(K)$ on $L$, respectively, and denote by $\vt=\vt_2\circ\vt_1$ their composition. Then $\eH^A=\eH^A_2\otimes\eH^A_1$ while $\eH^B=\eH^B_1\otimes\eH^B_2$. There is no possibility to discuss this away as, for instance, by arguments like $\vt_s\circ\vt_t=\vt_{s+t}=\vt_t\circ\vt_s$ when $\bfam{\vt_t}_{t\in\R_+}$ is an \nbd{E_0}semigroup. The corresponding isomorphism $\eH^A_2\otimes\eH^A_1\cong\eH^A_1\otimes\eH^A_2$ (or, similarly, for $\eH^B_i$) would not be the canonical one. A clear manifestation is Tsirelson's result \cite{Tsi00p1} that a product system of Hilbert spaces need not be isomorphic to its anti product system. (For Hilbert modules we may not even formulate what an anti product system is.)
\erem

\section{Generalizations of Bhat's approach}\label{Bha}

Under the hyposthesis of Theorem \ref{mainthm} (both for strict representations and for the \nbd{W^*}version) Bhat's approach generalizes easily, as shown in Skeide \cite{Ske02}, if $E$ has a \hl{unit vector} $\xi$, i.e.\ if $\AB{\xi,\xi}=\U$ (what, of course, includes that $\cB$ is unital). As in the proof for Hilbert spaces we define a Hilbert submodule $F_\xi=\vt(\xi\xi^*)F$ of $F$. As additional ingredient (as compared with Hilbert spaces) we define a left action of $\cB$ on $F_\xi$ by setting $by=\vt(\xi b\xi^*)y$. With these definitions one checks that
\beqn{
x\odot y
~\longmapsto~
\vt(x\xi^*)y
}\eeqn
defines an isometry $E\odot F_\xi\rightarrow F$. Like in the in proof of Theorem \ref{mainthm} surjectivity follows from existence of an approximate unit for $\sK(E)$ whose image under $\vt$ converges (strictly or \nbd{\sigma}weakly) to $\id_F$. Of course, $\vt(a)=a\odot\id_{F_\vt}$.

In this section we describe a construction from Skeide \cite{Ske03p1} which frees the preceding construction from the requirement of having a unit vector, at least, for the case of \nbd{W^*}modules. Then, as in Remark \ref{comprem}, we compare the construction with that one from Section \ref{intro}. Finally, we point out why the construction here, although not canonical (in the sense that it depends on the choice of a \it{complete quasi orthonormal system} for $E^*$), can have advantages over the intrinsic construction from Section \ref{intro}.

By making $\cB$ possibly smaller, we may always assure that $E$ is full and proofs of Theorem \ref{mainthm} which work for full $E$ work for arbitrary $E$. Existence of a unit vector is, however, a serious requirement. Our standard example is the \nbd{W^*}module $E=\Bfam{\substack{0\\\C\\\C}~\substack{\C\\0\\0}~\substack{\C\\0\\0}}\subset M_3$ which is a Hilbert module over $\cB=\bfam{\text{\raisebox{-.5ex}{$\substack{\smash{\C}\\\smash{0}}~\substack{\smash{0}\\\smash{M_2}}$}}}\subset M_3$ (with structures inherited from the embedding into $M_3$) which is full but does not admit a unit vector. (Actually, $E$ is a bimodule and as such $E$ is a Morita equivalence, because $\sB^a(E)=\sK(E)=\cB$.)

There are several equivalent possibilities to characterize \nbd{W^*}modules. A \nbd{W^*}module is always a Hilbert module $E$ over a \nbd{W^*}algebra $\cB$ fulfilling a further condition. We can require that $E$ be \it{self-dual} or that it has a predual Banach space. In the following section we consider von Neumann modules as introduced in Skeide \cite{Ske00b} as strongly closed operator spaces. A \nbd{W^*}module over a von Neumann algebra $\cB\subset\sB(G)$ is a von Neumann module and every von Neumann module is a \nbd{W^*}module. Many results on \nbd{W^*}modules have particularly simple and elementary proofs, when we transform them into von Neumann modules by choosing a faithful representation of $\cB$ on a Hilbert space $G$. Here we need the facts that $\sB^a(E)$ is a \nbd{W^*}algebra and that every \nbd{W^*}module admits a \hl{complete quasi orthonormal system}, i.e.\ a family $\bfam{e_\beta,p_\beta}_{\beta\in B}$ of pairs $(e_\beta,p_\beta)$ consisting of an element $e_\beta\in E$ and a projection $p_\beta\in\cB$ such that
\baln{
\AB{e_\beta,e_{\beta'}}
&
~=~
\delta_{\beta,\beta'}p_\beta
&
\text{and}
&&
\sum_{\beta\in B}e_\beta e_\beta^*
&
~=~
\id_E
}\ealn
where the sum is a \nbd{\sigma}weak limit over the increasing net of finite subsets of $B$. There is also a tensor product of \nbd{W^*}modules denoted by $\sbars{\odot}$.

So let us start with the assumptions of the \nbd{W^*}version of Theorem \ref{mainthm}. As explained before, we may assume that $E$ is full (which for \nbd{W^*}modules means the \nbd{\sigma}weakly closed ideal in $\cB$ generated by the range of the inner product of $E$ is $\cB$). It follows that the dual \nbd{\cB}\nbd{\sB^a(E)}module $E^*$ of $E$ is a Morita equivalence, in particular, that $\sB^a(E^*)=\cB$. Now choose a family $\bfam{e_\beta}_{\beta\in B}$ of elements in $E$ such that $\bfam{e_\beta^*,e_\beta e_\beta^*}_{\beta\in B}$ is a complete quasi orthonormal system for $E^*$. It follows that $p_\beta:=\AB{e_\beta,e_\beta}$ are projections in $\cB$ fulfilling $\sum_{\beta\in B}p_\beta=\U$.

\brem
If $B$ consists of a single element $\beta$, then $e_\beta$ is a unit vector. The following construction shows that the family $\bfam{e_\beta}_{\beta\in B}$, indeed, plays the role of the unit vector in the construction explain in the beginning of this section.
\erem

Now we define the \nbd{W^*}submodules $F_\beta=\vt(e_\beta e_\beta^*)F$ of $F$ and set $F_B=\bigoplus_{\beta\in B}F_\beta$. (The direct sum is that of \nbd{W^*}modules. Observe that the submodules $F_\beta$ of $F$ need not be orthogonal in $F$ so that $F_B$ is not a submodule of $F$.) On $F_B$ we define a left action of $b\in\cB$ by setting $by_\beta=\bigoplus_{\beta'\in B}\vt(e_{\beta'}be_\beta^*)y_\beta$ $(y_\beta\in F_\beta)$. (This defines, indeed, a \nbd{*}algebra representation of $\cB$ by adjointable operators on the algebraic direct sum, so the representing operators are bounded and, therefore, extend also to the \nbd{\sigma}weak closure.)

For $x\in E$ and $y_B\in F_B$ set $x_\beta=xp_\beta$ and $y_\beta=p_\beta y_B$. Then the mapping
\beqn{
x\odot y_B
~=~
\sum_{\beta\in B}x_\beta\odot y_\beta
~\longmapsto~
\sum_{\beta\in B}\vt(x_\beta e_\beta^*)y_\beta
}\eeqn
defines a unitary $E\sodots F_B\rightarrow F$ and $\vt(a)=a\odot\id_{F_B}$. The proof of isometry and surjectivity is exactly like in the version with a unit vector, except that now there is one index, $\beta$, more. See \cite{Ske03p1} for details.

\bemp[Comparison.]
In the case when there are unit vectors the comparison of $F_\xi$, $F_{\xi'}$ and $F_\vt$ works as in Remark \ref{comprem} ($\xi$ and $\xi'$ being possibily different unit vectors). $\vt(\xi'\xi^*)$ defines an isomorphism $F_\xi\rightarrow F_{\xi'}$ and $x^*\odot y\mapsto\vt(\xi x^*)y$ defines an isomorphism $F_\vt\rightarrow F_\xi$. For the family $\bfam{e_\beta}_{\beta\in B}$ the mapping
\beqn{
x^*\odot y
~=~
\sum_{\beta\in B}p_\beta x^*\odot y
~\longmapsto~
\bigoplus_{\beta\in B}\vt(e_\beta x^*)y
}\eeqn
defines an isomorphism $F_\vt\rightarrow F_B$. The identification of $F_B$ and $F_{B'}$ ($B'$ indicating the dual of some different complete quasi orthonormal system for $E^*$) follows by iterating the preceding formula with the inverse of its analogue for the other basis. The resulting formula is slightly complicated and does not give any new insight, so we do not write it down.
\eemp

\bemp[Advantages.]
If $\vt=\bfam{\vt_t}_{t\in\R_+}$ is an \hl{\nbd{E_0}semigroup} on $\sB^a(E)$ (i.e.\ a semigroup of unital strict or normal endomorphisms of $\sB^a(E)$), then the $E_t:=E_{\vt_t}=E^*\odot_t E$ form a product system in the sense of Bhat and Skeide \cite{BhSk00}. Indeed, if we identify $E_s\odot E_t$ with $E_{s+t}$ via
\beqn{
(x_s^*\odot_sy_s)\odot(x_t^*\odot_ty_t)
~\longmapsto~
x_s^*\odot_{s+t}\vt_t(y_sx_t^*)y_t,
}\eeqn
then
\baln{
(E_r\odot E_s)\odot E_t
&~=~
E_r\odot(E_s\odot E_t)
&
\text{and}
&&
(E\odot E_s)\odot E_t
&~=~
E\odot(E_s\odot E_t).
}\ealn
Also the identifications via a unit vector $\xi$ or a family $\bfam{e_\beta}_{\beta\in B}$ respect these associativity conditions. So far, the two constructions can be used interchangeably. This changes, however, when we wish to include also technical conditions on product systems.

A product system of Hilbert spaces in the sense of Arveson \cite{Arv89} is supposed to be derived from an \nbd{E_0}semigroup on $\sB(H)$ that is pointwise \nbd{\sigma}weakly continuous (in time). The product system has, therefore, the structure of a Banach bundle, more precisely, the structure of a trivial Banach bundle.

In Skeide \cite{Ske03b} we have investigated the Hilbert module version in presence of a unit vector. Requiring the product system to be isomorphic to a trivial Banach bundle seems too much. (We do not even know, whether all members $E_t$ ($t>0$) of a product system are isomorphic as right modules.) However, our product systems have the structure of a subbundle of a trivial Banach bundle. This can be derived easily from the observation that in presence of a unit vector all $E_t$ can be identified with submodules $\vt_t(\xi\xi^*)E$ of $E$. Since $\vt$ is sufficiently continuous, the corresponding subbundle of the trivial Banach bundle $\RO{0,\infty}\times E$ is a Banach bundle (there are enough \it{continuous sections}).

A \nbd{\sigma}weak version in presence of a unit vector does not seem to present a difficulty. Now, if we have a family $\bfam{e_\beta}_{\beta\in B}$ instead of a unit vector, $F_B$ need no longer be a submodule of $E$. However, each $F_\beta=\vt(e_\beta e_\beta^*)F$ is a submodule. It follows that $F_B=\bigoplus_{\beta\in B}F_\beta$ is a submodule of $\bigoplus_{\beta\in B}F$. Therefore, it seems reasonable to expect that the product system is a \nbd{\sigma}weak subbundle of the trivial \nbd{\sigma}weak bundle $\RO{0,\infty}\times\bigoplus_{\beta\in B}E$. This requires a convenient definition of \nbd{\sigma}weak bundle and is work in progress.

In both cases the construction according to Section \ref{intro} does not seem to help to identify a good candidate for the trivial bundle of which the product system is a subbundle.
\eemp

\section{The generalization of Arveson's approach}\label{Arv}

For this section we need a longer preparation. If $E$ is a \nbd{\cB}\nbd{\cB}module, then the \hl{\nbd{\cB}center} of $E$ is the space
\beqn{
C_\cB(E)
~=~
\CB{x\in E\colon bx=xb~(b\in\cB)}.
}\eeqn
In what follows it is essential that von Neumann algebras and von Neumann modules (or, more generally, Hilbert modules over von Neumann algebras) always come along with an identification as concrete subspaces of operators on or between Hilbert spaces. A von Neumann algebra $\cB$ is given as a concrete subalgebra of $\sB(G)$ acting (always nondegenerately) on a Hilbert space $G$. Every Hilbert \nbd{\cB}module $E$ may, then, be identified as a \nbd{\cB}submodule of $\sB(G,H)$ for a suitable Hilbert space $H$ in the following way.

Set $H=E\odot G$. Then, an element $x\in E$ defines an operator $L_x\colon g\mapsto x\odot g$ in $\sB(G,H)$. Clearly, $\AB{x,y}=L_x^*L_y$ and $L_{xb}=L_xb$. Moreover, $E$ acts nondegenerately on $G$ in the sense that $L_EG$ is total in $H$ and the pair $H,\eta\colon x\rightarrow L_x$ is determined by these properties up to (unique) canonical isomorphism. We, therefore, identify $E$ as a subset of $\sB(G,H)$ by identifying $x$ with $L_x$. Following Skeide \cite{Ske00b} $E$ is a \hl{von Neumann module}, if it is strongly closed in $\sB(G,H)$. One may show (see \cite{Ske00b}) that a Hilbert \nbd{\cB}module $E$ over a von Neumann algebra $\cB\subset\sB(G)$ is self-dual, if and only if $E$ is a von Neumann module.

On $H$ we have a normal unital representation $\rho'$ of $\cB'$, the \hl{commutant lifting}, defined by $\rho'(b')=\id_E\odot b'$. The space $C_{\cB'}(\sB(G,H))$ is a von Neumann \nbd{\cB}module containing $E$ as a submodule with zero-complement. Since $C_{\cB'}(\sB(G,H))$ is self-dual, it follows that $E=C_{\cB'}(\sB(G,H))$, if and only if $E$ is a von Neumann module. Observe that in this case $\rho'(\cB')'$ is exactly $\sB^a(E)$.

The identification of $C_{\cB'}(\sB(G,H))$ as the unique minimal self-dual extension of $E$ (in the sense of Paschke \cite{Pas73}) was already known to Rieffel \cite{Rie74a}. The definition of von Neumann modules seems to be due to \cite{Ske00b}. In Skeide \cite{Ske03p} we show directly (without self-duality of von Neumann modules) that $E=C_{\cB'}(\sB(G,H))$, if $E$ is a von Neumann module, and then give a different proof of Rieffel's result that $C_{\cB'}(\sB(G,H))$ is self-dual. Muhly and Solel \cite{MuSo02} show that, conversely, every normal unital representation $\rho'$ of $\cB'$ on a Hilbert space gives rise to a von Neumann \nbd{\cB}module $C_{\cB'}(\sB(G,H))\subset\sB(G,H)$ acting nondegenerately on $G$. Summarizing, we have a one-to-one correspondence (up to canonical isomorphisms)
\beqn{
\sB(G,H)
~\supset~
E
~~~\longleftrightarrow~~~
(\rho',H)
}\eeqn
between von Neumann \nbd{\cB}modules and normal representations of $\cB'$.

If $E$ is a von Neumann \nbd{\cA}\nbd{\cB}module (that is, $\cA\subset\sB(K)$ is another von Neumann algebra and the canonical homomorphism $\cA\rightarrow\sB^a(E)\rightarrow\sB(H)$ defines a normal unital representation $\rho$ of $\cA$ on $H$), then we have a pair of representations $\rho$ and $\rho'$ with mutually commuting ranges. $\rho'$ gives back the right module $E$ as intertwiner space $C_{\cB'}(\sB(G,H))$ and $\rho$ gives back the correct left action. This works also if we start with a triple $(\rho,\rho',H)$. For the standard representation of $\cB$ so that $\cB'\cong\cB^{op}$ and $\rho'$ may be viewed as representation of $\cB^{op}$, we are in the framework of Connes and others where von Neumann bimodules and pairs of representations of $\cA$ and $\cB^{op}$ are interchangeable pictures of the same thing. The more general setting where $\cB$ is not necessarily given in standard representation seems not to have been observed before Skeide \cite{Ske03c} and Muhly and Solel \cite{MuSo03p}. Going only slightly further, by exchanging the roles of $\cA$ and $\cB$ in the triple $(\rho,\rho',H)$, we find the following one-to-one correpondence
\beqn{
\parbox{8cm}{
\xymatrix{
	&(\rho,\rho',H)	\ar@{<->}[dl]	\ar@3{-}[r]	&(\rho',\rho,H)	\ar@{<->}[dr]	&	\\
E	&									&									&E'
}
}
}\eeqn
where $E'=C_\cA(\sB(K,H))$ is a von Neumann \nbd{\cB'}\nbd{\cA'}module. We refer to $E'$ as the \hl{commutant} of $E$ (and conversely), because when $E$ is the von Neumann \nbd{\cB}\nbd{\cB}module $\cB$, then $E'=\cB'$. Also this correspondence was observed in Skeide \cite{Ske03c} and, later, in Muhly and Solel \cite{MuSo03p}. See also Gohm and Skeide \cite{GoSk03p} for another application of the commutant.

\bob\label{impob}
It is important to notice that the preceding correpondences between (bi-)mod\-ules and (pairs of) representations enables us to identify von Neumann (bi-)modules by, first, identifying Hilbert spaces and, then, showing that representations on them coincide.
\eob

\bob
$E=C_{\cB'}(\sB(G,H))$ and $E'=C_\cA(\sB(K,H)$ act nondegenerately on $G$ and $K$, respectively. Therefore, $\cls EG=H=\cls E'K$. Since we canonically identify $H=E\odot G$ and $H=E'\odot K$ (by setting $xg=x\odot g$ and $x'k=x'\odot k$), we have $E\odot G=E'\odot K$. Writing down this identity is an invitation to the reader to take an element $x\odot g$ in $E\odot G$ and write it as a sum of elements $x'\odot k$ in $E'\odot K$. There is no canonical way how to do it, like there is no canonical way how to express a general element in a tensor product by a sum over elementary tensors. We just know that it is possible and that how ever we do it our conclusions do ot depend on the choice.

For instance, it is important to keep in mind how the representations $\rho$ and $\rho'$ act in these pictures. We have $\rho(a)(x\odot g)=ax\odot g$, while $\rho(a)(x'\odot k)=x'\odot ak$ and, conversely, $\rho'(b')(x\odot g)=x\odot b'g$, while $\rho'(b')(x'\odot k)=b'x'\odot k$.
\eob

Now we come to the third proof of Theorem \ref{mainthm} where we, acutally, first construct the commutant of $F_\vt$. We assume the hypothesis for the \nbd{W^*}version of Theorem \ref{mainthm}. As in Section \ref{Bha} we assume that $E$ is full. Furthermore, we assume that $\cB\subset\sB(G)$ and $\cC\subset\sB(L)$ so that $E$ and $F$ are von Neumann modules. We make up the following dictionary.
\baln{
H		~&=~E\odot G		&	K		~&=~F\odot L		&&\\
\rho'(b')	~&=~\id_E\odot b'		&	\sigma'(c')	~&=~\id_F\odot c'		&&(b'\in\cB',c'\in\cC')\\
\rho(a)	~&=~a\odot\id_G		&	\sigma(a)	~&=~\vt(a)\odot\id_L	&&(a\in\sB^a(E))
}\ealn
It makes, therefore, sense to define the intertwiner space $F'_\vt=C_{\sB^a(E)}(\sB(H,K))$ which is the subspace of $\sB(G,K)$ of all mappings intertwining the actions of $\sB^a(E)$ via $\sigma$ and $\rho=\id_{\sB^a(E)}$ where by definition $\sB^a(E)$ is a von Neumann algebra on $H$ via the identity representation $\rho$.

Recall that the commutant of $\sB^a(E)$ is $\rho'(\cB')$. Therefore, in the above correpondence between von Neumann bimodules and pairs of representations, we may consider $F'_\vt$ as the von Neumann \nbd{\cC'}\nbd{\rho'(\cB')}module determined by the triple $(\sigma',\sigma,K)$ with inner product ${y'_1}^*y'_2\in\rho'(\cB')$, and left and right multiplication given simply by composition with $\sigma'(c')$ from the left and with $\rho'(b')$ from the right, respectively.

Now, since $E$ is full, $\rho'$ is faithful so that $\rho'(\cB')\cong\cB'$. Therefore, we may, finally, and will consider $F'_\vt$ as von Neumann \nbd{\cC'}\nbd{\cB'}module where
\beqn{
\AB{y'_1,y'_2}
~:=~
{\rho'}^{-1}({y'_1}^*y'_2)
\text{~~~~~~and~~~~~~}
c'y'b'
~:=~
\sigma'(c')y'\rho'(b').
}\eeqn

The following identification
\beqn{
F\odot L
~=~
K
~=~
F'_\vt\odot H
~=~
F'_\vt\odot E\odot G
~=~
E\odot F'_\vt\odot G
~=~
E\odot F''_\vt\odot L
}\eeqn
identifies the Hilbert spaces $F\odot L$ and $E\odot F''_\vt\odot L$ of the von Neumann modules $F$ and $E\odot F''_\vt$.

\bob\label{EE'}
The ``tricky'' identification is $F'_\vt\odot E\odot G=E\odot F'_\vt\odot G$. One easily checks that there is a canonical identification of these spaces simply by flipping the first two factors in elementary tensors. This also shows that operators on $E$ and on $F'_\vt$, respectively, act directly on the factor where they belong. In \cite{Ske03p1} we investigate systematically the tensor product $E\sodots E'\cong E'\sodots E$ of a von Neumann \nbd{\cB}module and a von Neumann \nbd{\cB'}module which is a von Neumann \nbd{(\cB\cap\cB')'}module. This tensor product may be viewed as a generalization of the exterior tensor product with which it has many properties in common.
\eob

An investigation how the relevant algebras act on these Hilbert spaces show that the von Neumann \nbd{\cC}modules $F$ and $E\sbars{\odot}F''_\vt$ coincide (in the sense of Observation \ref{impob}) and that $\vt(a)=a\odot\id_{F_\vt}$ what concludes the third proof.

\bemp[Applications.]
Taking the commutant of von Neumann bimodules is anti-multiplicative under tensor product; see \cite{Ske03p1} for details. Taking into account that in Section \ref{B(G)}, clearly, $\eH^A=(\eH^B)'$, we see that the Arveson system of an \nbd{E_0}semigroup on $\sB(H)$ is the opposite of its Bhat system. Also the product systems of \nbd{\cB}\nbd{\cB}modules in \cite{BhSk00} and of \nbd{\cB'}\nbd{\cB'}modules in \cite{MuSo02}, both constructed from the same CP-semigroup on $\cB$, are commutants of each other. We explain this in \cite{Ske03c}.

Also other applications are related to endomorphisms of $\sB^a(E)$. While every bimodule $F_\vt$ comes from a representation of $\sB^a(E)$ on $E\odot F_\vt$, the question, whether a bimodule comes from an endomorphism (i.e.\ whether there exists an $E$ such that $E\odot F_\vt\cong E$) is nontrivial. It is equivalent to the question whether $F'_\vt$ has an isometric fully coisometric covariant representation on a Hilbert space. In the semigroup version this means that the question, whether a product system stems from an \nbd{E_0}semigroup on some $\sB^a(E)$, is equivalent to the question, whether the commutant system allows for such a covariant representation.

We investigate these and other questions in Muhly, Skeide and Solel \cite{MSS02p}.
\eemp

\bemp[Comparison.]
How is $F''_\vt$ related to $F_\vt=E^*\sodots F$ from Section \ref{intro}? Of course, we know that they are canonically isomorphic, but we want to see the identification in the sense of Observation \ref{impob}. In fact, we are able to identify $(E^*\odot F)'=F'\odot{E^*}'$ and $F'_\vt$, but after the sketchy discussion earlier in this section it is not possible to present the subtle arguments (flipping continuously between the isomorphic von Neumann algebras $\cB'$ and $\rho'(\cB')$) in a coherent way. (In fact, many readers will feel uncomfortable with our continuously used canonical identifications of spaces which \it{a priori} are different, and doing this consistently requires a skillful preparation.) Once more, we refer the reader to \cite{Ske03p1} for a detailed discussion.
\eemp

\lf\noindent
\bf{Acknowledgements.} The results of Section \ref{intro} are joint work with Paul Muhly and Baruch Solel and most of these and other results have been worked out during the author's stays at ISI Bangalore and University of Iowa in 2003. The author wishes to express his gratitude for hospitality during two fantastic stays to B.V.Rajarama Bhat (ISI) and Paul S.\ Muhly (University of Iowa).

\setlength{\baselineskip}{2.5ex}


\newcommand{\Swap}[2]{#2#1}\newcommand{\Sort}[1]{}
\providecommand{\bysame}{\leavevmode\hbox to3em{\hrulefill}\thinspace}
\providecommand{\MR}{\relax\ifhmode\unskip\space\fi MR }
\providecommand{\MRhref}[2]{%
  \href{http://www.ams.org/mathscinet-getitem?mr=#1}{#2}
}
\providecommand{\href}[2]{#2}


\end{document}